\documentclass[11pt,a4paper]{amsart}
\usepackage[all]{xy}
\usepackage{amssymb}
\usepackage{amsmath}
\usepackage{bm}
\usepackage[colorlinks,linkcolor=blue]{hyperref}
\usepackage{color,xcolor,dsfont}
\usepackage{enumerate}
\usepackage{epsfig}
\usepackage{extarrows}
\usepackage{framed}
\usepackage{graphicx}
\usepackage{marvosym}
\usepackage{mathrsfs}
\usepackage{setspace} 
\usepackage{tikz}
\usepackage{tikz-cd}

\newcommand{\scha}{\mathcal{A}}

\newcommand{\schd}{\mathcal{D}}

\newcommand{\she}{\mathscr{E}}
\newcommand{\shf}{\mathscr{F}}

\newcommand{\shh}{\mathscr{H}}
\newcommand{\shi}{\mathscr{I}}
\newcommand{\shk}{\mathscr{K}}

\newcommand{\sho}{\mathscr{O}}
\newcommand{\shp}{\mathscr{P}}

\newcommand{\shu}{\mathscr{U}}

\newcommand{\mukai}{\textup{\textsf{v}}}
\newcommand{\deriver}{\textup{\textsf{R}}}
\newcommand{\derivel}{\textup{\textsf{L}}}

\DeclareMathOperator{\ev}{ev}

\DeclareMathOperator{\Aut}{Aut}

\DeclareMathOperator{\Hom}{Hom}

\DeclareMathOperator{\image}{im}
\DeclareMathOperator{\spec}{Spec}

\newcommand{\bz}{\mathbb{Z}}
\newcommand{\bq}{\mathbb{Q}}

\newcommand{\bc}{\mathbb{C}}

\newcommand{\bp}{\mathbb{P}}

\newcommand{\ch}{\textup{\textsf{ch}}}
\newcommand{\td}{\textup{\textsf{td}}}

\newcommand{\defi}[1]{{\textbf{\emph{#1}}}}
\newcommand{\cate}[1]{{\textbf{\textup{#1}}}}
\newcommand{\funct}[1]{{\textup{\textbf{\textsf{#1}}}}}

\theoremstyle{plain}
\newtheorem{theorem}{Theorem}[section]
\newtheorem{corollary}[theorem]{Corollary}
\newtheorem{lemma}[theorem]{Lemma}
\newtheorem{proposition}[theorem]{Proposition}

\theoremstyle{definition}
\newtheorem{definition}[theorem]{Definition}

\newtheorem{example}[theorem]{Example}

\author{Ziqi Liu}
\title{On Fourier--Mukai type autoequivalences of Kuznetsov components of cubic threefolds}
\address{Mathematisches Institut, Universitat Bonn, Endenicher Allee 60, 53115 Bonn, Germany}
\email{ s6ziliuu@uni-bonn.de}

\begin{document}

\maketitle

\vspace{-8mm}

\begin{abstract}
We determine the group of all Fourier--Mukai type autoequivalences of Kuznetsov components of smooth complex cubic threefolds, and provide yet another proof for the Fourier--Mukai version of categorical Torelli theorem for smooth complex cubic threefolds.
\end{abstract}

\section{Introduction}
\subsection{}
The seminal results of Bondal and Orlov motivate the study on derived categories of smooth projective varieties. This interest later spreads to other triangulated categories, and an important viewpoint is that a nice enough triangulated category can be regarded as (the derived category of) a noncommutative variety.

Among these triangulated categories, the Kuznetsov component $\scha_Y$ of a (smooth complex) cubic threefold $Y$ is one of the most studied ones. It carries enough information to determine $Y$ (see, for example, \cite{BBFHMRS23,BMMS12,FLZ23,PY22}). However, many geometric aspects for $\scha_Y$ as a noncommutative variety have not been well-understood.

In this paper, we study the autoequivalences of $\scha_Y$, and particularly determine all Fourier--Mukai type ones. Also, we give an elementary proof for the Fourier--Mukai version of categorical Torelli theorem for cubic threefolds.

\subsection{}
Given a cubic threefold $Y$, we denote by $\Aut_{FM}(\scha_Y)$ the set of all Fourier--Mukai type autoequivalences of $\scha_Y$. According to \cite{HR19}, many natural autoequivalences, such as the identity functor $\funct{Id}_{\scha_Y}$ and the Serre functor $\funct{S}$, are in $\Aut_{FM}(\scha_Y)$. So one can see that $\Aut_{FM}(\scha_Y)$ is really a group under the functor composition.

Set $j^*\colon D^b(Y)\rightarrow\scha_Y$ to be the left adjoint functor of the inclusion $j_*\colon\scha_Y\subset D^b(Y)$, then the image $j^*k(x)$ is isomorphic to a stable sheaf $\shk^x$ on $Y$ up to a shift. Thanks to \cite{BBFHMRS23,BMMS12}, we are able to make a decomposition
$$\Aut_{FM}(\scha_Y)\cong\bz\times\Aut(Y)\times G_0$$
where $G_0:=\{\Phi\in\Aut_{FM}(\scha_Y)\,|\,\Phi(\shk^x)\cong \shk^x\,\,\textup{ for each point }x\in Y\}$. To decode $G_0$, we need to understand the relationship between certain flat families of stable sheaves on $Y$ with class $v=v(\shk^x)$. It turns out that $G_0$ is trivial, due to the observation that the moduli space $M_Y(v)$ of stable sheaves considered in \cite{BBFHMRS23} is indeed a fine moduli space and admits a universal sheaf.

\begin{theorem}[Theorem \ref{5_main-theorem}]\label{0_main-theorem}
One has $\Aut_{FM}(\scha_Y)\cong\bz\times\Aut(Y)$.
\end{theorem}

This result is obtained independently also in \cite{FLZ23} among other things, using a complete different method. Explicitly, they prove that each element in $\Aut_{FM}(\scha_Y)$ can be extended to an autoequivalence of $D^b(Y)$ and then conclude by \cite[Proposition 4.17]{HuFu}.

\subsection{}
It is believed that any equivalence between Kuznetsov components of cubic threefolds is of Fourier--Mukai type, which is indeed the case for cubic fourfolds cf.\ \cite{LPZ23}. Therefore, it should be sufficient to consider the following form of categorical Torelli theorem for cubic threefolds.

\begin{theorem}[Theorem \ref{5_another-proof}]\label{0_another-proof}
Consider two cubic threefolds $Y_1$ and $Y_2$, and assume that there exists a Fourier--Mukai type equivalence $\scha_{Y_1}\cong\scha_{Y_2}$, then $Y_1\cong Y_2$.
\end{theorem}

This specific form of categorical Torelli theorem has been proved in \cite{BT16,Pe22} using certain motivic constructions of intermediate Jacobians, and proved in \cite{LZ23} using matrix factorization. Here we provide another proof for this theorem based on \cite{OR18} and the spirit of \cite{HR19}. Our argument is similar to that for \cite[Proposition 5.23]{Pe22} but more elementary.

\subsection{Organization} In Section \ref{SS_Kuznetsov components}, we review definitions, clear necessary notations and exposit a couple of recent results about Kuznetsov components of cubic threefolds. In Section \ref{SS_Fourier--Mukai}, we recall some aspects of Fourier--Mukai transforms. Section \ref{SS_Main Theorem} is devoted to Theorem \ref{0_main-theorem} and we give our proof for Theorem \ref{0_another-proof} in Section \ref{SS_FM categorical Torelli}.

\subsection{Conventions} 
In this paper, all varieties are defined over $\bc$, all triangulated categories are assumed to be complex linear, and functors between triangulated categories are asked to be over $\bc$. Given a variety $X$, we set $D^b(X):=D^b(\cate{Coh}(X))$ to be the bounded derived category of coherent sheaves on $X$. By a cubic threefold, we always mean a smooth projective one, which will be denoted by $Y$ or $Y_i$; by a point, we always mean a closed point.

\subsection*{Acknowledgements}
This is a modified version of the first part of my Master's thesis at the University of Bonn. I would like to thank my supervisor Daniel Huybrechts for his patience and guidance, and thank Shizhuo Zhang for many helpful discussions and comments. I also wish to thank the referees for careful reading of the paper and for pointing out typos and inaccuracies.

\section{Kuznetsov components of cubic threefolds}\label{SS_Kuznetsov components}
Here we introduce some basic aspects of Kuznetsov components of cubic threefolds and a couple of important results. Standard notions and notations can be found in \cite{HuFu,HuCu}.

\subsection{}
Let $\iota\colon Y\subset\bp^4$ be a cubic threefold polarized by $\sho_Y(1):=\iota^*\sho_{\bp^4}(1)$. The \emph{Kuznetsov component} $\scha_Y$ of $Y$ is defined to be the right orthogonal complement
$$\scha_Y:=\{E\in D^b(Y)\,|\,\Hom(\sho_Y,E[n])=\Hom(\sho_Y(1),E[n])=0\textup{ for all }n\in\bz\}$$
of the exceptional collection $\sho_Y,\sho_Y(1)$ in $D^b(Y)$.

\begin{example}
Let $\ell\subset Y$ be any line, then the ideal sheaf $\shi_{\ell}$ lies in $\scha_Y$ by definition.
\end{example}

The subcategory $\scha_Y\subset D^b(Y)$ is admissible, and we will denote the inclusion functor and its right (resp.\ left) adjoint functor by $j_*$ and $j^!$ (resp.\ $j^*$). Explicitly, one has
$$j^*\colon D^b(Y)\stackrel{\funct{L}_{\sho_Y(1)}}{\longrightarrow}\langle\scha_Y,\sho_Y\rangle\stackrel{\funct{L}_{\sho_Y}}{\longrightarrow}\scha_Y$$
where $\funct{L}_{E_0}\colon \schd\rightarrow \langle E_0\rangle^{\perp}$ is the left mutation functor represented by the mapping
$$E\mapsto\textsf{C}(\ev\colon \bigoplus_m\Hom(E_0,E[m])\otimes E_0[-m]\rightarrow E)$$
for any exceptional object $E_0\in\schd$ of any given triangulated category.

\begin{example}\label{2_projection-functor}
Let $k(x)$ be the skyscraper sheaf at a point $x\in Y$, then 
$$j^*k(x)=\funct{L}_{\sho_Y}\circ\funct{L}_{\sho_Y(1)}k(x)=\funct{L}_{\sho_Y}\shi_x(1)[1]=\textup{\textsf{C}}\left(\sho_Y\otimes H^0(Y,\shi_{x}(1))\rightarrow\shi_{x}(1)\right)[1]\cong\shk^x[2]$$
where $\shi_x$ is the ideal sheaf of the point $x$ and $\shk^x$ is the coherent sheaf defined by the short exact sequence $0\rightarrow\shk^x\rightarrow\sho_Y\otimes H^0(Y,\shi_x(1))\rightarrow\shi_x(1)\rightarrow0$.
\end{example}

\subsection{}
The Kuznetsov component $\scha_Y$ has many natural non-trivial autoequivalences.

\begin{example}
Let $f\colon Y\rightarrow Y$ be an automorphism, then the derived direct image $\deriver f_*$ and the derived inverse image $\derivel f^*$ restrict to mutually inverse autoequivalences of $\scha_Y$.
\end{example}

\begin{example}
Since $\scha_Y\subset D^b(Y)$ is admissible, it has a Serre functor $\funct{S}$ which can be represented by $E\mapsto j^!(j_*E\otimes\sho_Y(-2))[3]$. Accordingly, one has $\funct{S}^{-1}\cong j^*(j_*(-)\otimes\sho_Y(2))[-3]$.
\end{example}

\begin{example}\label{2-2_degree-shift}
The tensor product with a line bundle on $Y$ does not restrict to an autoequivalence of $\scha_Y$. Instead, one has the \emph{degree shift} functor $\funct{T}\colon E\mapsto j^*(j_*E\otimes\sho_Y(1))$ such that
$$\funct{T}^n\cong j^*(j_*(-)\otimes\sho_Y(n))$$
for $n\geq0$. In particular, one has $\funct{S}^{-1}\cong\funct{T}^{2}\circ[-3]$. 
\end{example}

\subsection{} It is proved in \cite{Ku04} that $\funct{T}^3\cong[2]$, hence $\funct{S}\cong\funct{T}\circ[1]$ and $\funct{S}^3\cong [5]$. This property allows one to give an explicit description of the numerical Grothendieck group $N(\scha_Y)$.
 
\begin{proposition}[{\cite[Proposition 2.7]{BMMS12}}]\label{2_numerical-Grothendieck-group}
Consider the ideal sheaf $\shi_{\ell}$ of a line $\ell\subset Y$, then 
$$N(\scha_Y)\cong\bz[\shi_{\ell}]\oplus\bz[\cate{S}\shi_{\ell}]$$ 
and the Euler pairing $\chi(-,-)$ on $N(\scha_Y)$ has the following form 
$$\begin{pmatrix}
		-1&-1\\0&-1
\end{pmatrix}$$
with respect to this basis.
\end{proposition}

\begin{example}
One computes that $[\funct{S}^2\shi_{\ell}]=-[\shi_{\ell}]+[\funct{S}\shi_{\ell}]$ and $[j^*k(x)]=[\shi_{\ell}]+[\funct{S}\shi_{\ell}]$.
\end{example}

Also in \cite{BMMS12}, they construct the first (Bridgeland) stability condition on $\scha_Y$ and prove that the Kuznetsov component $\scha_Y$ is an invariant for cubic threefolds. Here we want to mention the following proposition, which will be used later.

\begin{proposition}\label{2-3_decomposition}
An autoequivalence $\Phi\colon\scha_Y\cong\scha_{Y}$ is decomposed into the form $\funct{S}^r\circ[m]\circ\Phi_0$ for some integers $r,m\in\bz$ and an autoequivalence $\Phi_0\colon\scha_Y\cong\scha_{Y}$ which sends the ideal sheaf $\shi_{\ell}$ of a line on $Y$ to the ideal sheaf $\shi_{\ell'}$ of some line $\ell'$ on $Y$.

\begin{proof}
This follows from the proof for \cite[Proposition 5.1]{BMMS12}.
\end{proof}
\end{proposition}

\subsection{} 
Given a (numerical) stability condition $\sigma$ on $\scha_Y$, people are interested in the moduli space $M_{\sigma}(v)$ of $\sigma$-semistable objects with numerical class $v\in N(\scha_Y )$. Concerning the Kuznetsov components of cubic threefolds, the first such moduli space is essentially constructed in \cite{BMMS12}: the moduli space $M_{\sigma}([\shi_{\ell}])$ is isomorphic to the Fano variety of lines on $Y$ for any Serre-invariant stability condition on $\scha_Y$ cf.\ \cite{PY22}.

Bayer et al \cite{BBFHMRS23} study another family of objects and obtain that

\begin{theorem}
Given any Serre-invariant stability condition $\sigma$ on $\scha_Y$, then one has
$$M_{\sigma}([\shk^x])=\{\shk^x\,|\,x\in Y\}\cup\{\she_C\,|\,C\subset Y\textup{ is a twisted cubic}\}$$
as a set. One also notices that all the $\shk^x$'s and $\she_C$'s are $\sigma$-stable.\footnote{Here $\she_C$ is a stable locally free sheaf on $Y$ lying in $\scha_Y$. We omit its construction since we will not use it.}
\end{theorem}

Explicitly, they show that $M_{\sigma}([\shk^x])$ is identified with the moduli space $M_Y(v)$ of Gieseker semistable sheaves on $Y$ with Chern character $v=\ch(\shk^x)$. Moreover, the sheaves in $M_Y(v)$ are all stable as $\ch(\shk^x)=(3,-h,-\frac{1}{2}h^2,\frac{1}{6}h^3)$ for $h:=c_1(\sho_Y(1))$.

Using a certain universal sheaf $\shu$ on $Y\times Y$, they construct a closed embedding $i\colon Y\hookrightarrow M_Y(v)$. Under this embedding, they prove that an automorphism $M_{Y}(v)\cong M_{Y}(v)$ descends to an automorphism $Y\cong Y$ cf.\ \cite[Lemma 7.5]{BBFHMRS23}. Combined with $M_{\sigma}([\shk^x])\cong M_Y(v)$, one sees

\begin{proposition}\label{2-3_point-identification}
An autoequivalence $\Phi\colon\scha_Y\cong\scha_Y$ induces a bijective map 
$$\{\shk^x\,|\,x\in Y\}\rightarrow\{\funct{S}^r\shk^{x}\,|\,x\in Y\}$$
for some fixed $r$ depending only on $\Phi$.
\end{proposition}

To conclude this section, we make the observation that $M_Y(v)$ is a fine moduli space. Indeed, one can use Hirzebruch--Riemann--Roch formula to compute that
\begin{align*}
\chi(\shk^x\otimes\sho_{\ell})=&\int_Y(3,-h,-\frac{1}{2}h^2,\frac{1}{6}h^3)\cup(0,0,\frac{1}{3}h^2,0)\cup(1,h,\frac{2}{3}h^2,\frac{1}{3}h^3)=2\\
\chi(\shk^x\otimes\sho_Y(2))=&\int_Y(3,-h,-\frac{1}{2}h^2,\frac{1}{6}h^3)\cup(1,2h,2h^2,\frac{4}{3}h^3)\cup(1,h,\frac{2}{3}h^2,\frac{1}{3}h^3)=27
\end{align*}
and then concludes by \cite[Theorem 4.6.5]{HL}.

\section{Fourier--Mukai transforms}\label{SS_Fourier--Mukai}
Here we introduce some necessary facts about Fourier--Mukai transforms. We will tactically use the analytic topology when we talk about cohomology and topological $K$-theory and will not notationally distinguish a variety and its analytification.

\subsection{}
At first we collect some standard materials from \cite{HuFu}. Let $X$ and $Z$ be two smooth projective varieties and choose an object $P\in D^b(X\times Z)$, we define
$$\Phi_P\colon D^b(X)\rightarrow D^b(Z),\quad E\mapsto \deriver p_*(P\otimes^{\derivel} \derivel q^*E)$$
to be the \emph{Fourier--Mukai transform with kernel $P$} where $p\colon X\times Z\rightarrow Z$ and $q\colon X\times Z\rightarrow X$ are the projections. Here $\derivel q^*$ is indeed just the non-derived $q^*$ as the projection $q$ is flat.

The composition $\Phi_P\circ\Phi_Q$ of two Fourier--Mukai transforms $\Phi_Q\colon D^b(X)\rightarrow D^b(Z)$ and $\Phi_P\colon D^b(Z)\rightarrow D^b(W)$ is isomorphic to a Fourier--Mukai transform $\Phi_{P\circ Q}$ with kernel
$$P\circ Q:=\deriver\pi_{XW*}(\pi^*_{ZW}P\otimes^{\derivel}\pi^*_{XZ}Q)$$ 
where $\pi_{XZ},\pi_{ZW},\pi_{XW}$ are the projections from $X\times Z\times W$ onto the respective factors.

The Fourier--Mukai transform $\Phi_P\colon D^b(X)\rightarrow D^b(Z)$ admits a left adjoint functor $\Phi_{P_\textrm{L}}$ and a right adjoint functor $\Phi_{P_{\textrm{R}}}$ with kernels
$$P_{\textrm{L}}\cong P^\vee\otimes p^*\omega_Z[\dim{Z}]\quad\textup{and}\quad P_{\textrm{R}}\cong P^\vee\otimes q^*\omega_X[\dim{X}]$$
where $P^\vee:=\deriver\shh om(P,\sho_{X\times X})$ is the derived dual.

\subsection{}
The idea of Fourier--Mukai transforms can be applied to maps between other invariants in complex algebraic geometry. The most commonly used ones are cohomological Fourier--Mukai maps. Under the same assumption as before, one defines
$$\Phi^H_{P}\colon H^{\bullet}(X,\bq)\rightarrow H^{\bullet}(Z,\bq),\quad v\mapsto p_*(q^*v\cup \mukai(P))$$
where $p_*,q^*$ are pushforward and pullback in cohomology and $\mukai(P):=\ch([P])\cup\sqrt{\td(X\times Z)}$ is the Mukai vector of the kernel $P\in D^b(X\times Z)$. Though cohomological Fourier--Mukai maps do not preserve graded structure on $H^{\bullet}$ in general, one has
$$\Phi_P^H( H^{p,q}(X))\subset\bigoplus_{p-q=r-s} H^{r,s}(Z)$$
In particular, a group isomorphism $\Phi_P^H\colon H^{\bullet}(X,\bq)\rightarrow H^{\bullet}(Z,\bq)$ will restrict to isomorphisms
$$\bigoplus_{p-q=i} H^{p,q}(X)\cong\bigoplus_{p-q=i} H^{p,q}(Z)$$
Moreover, one always has a commutative diagram
\begin{displaymath}\xymatrix{
		D^b(X)\ar[d]_{\Phi_P}\ar[rr]^{\mukai}&&H^{\textup{even}}(X,\bq)\ar[d]_{\Phi_P^H}\\
		D^b(Z)\ar[rr]^{\mukai}&&H^{\textup{even}}(Z,\bq)
}\end{displaymath}

\subsection{}
Similar to the cohomological Fourier--Mukai maps, one can also consider the topological $K$-theoretic Fourier--Mukai maps cf.\ \cite{OR18}. Under the same assumption as before, one defines
$$\Phi^K_{P}\colon K^0(X)\rightarrow K^0(Z),\quad \kappa\mapsto p_*(q^*\kappa\cup [P])$$
where $p_*,q^*$ are pushforward and pullback in topological $K$-theory and $[P]\in K^0(X\times Z)$ is the image of $P$ under the natural map $D^b(X\times Z)\rightarrow K^0_{a}(X\times Z)\rightarrow K^0(X\times Z)$.

The pullback in topological $K$-theory is induced by the pullback of vector bundles, while the pushforward is not canonical. It is constructed by Atiyah and Hirzebruch in \cite{AH61} and elaborated on by Karoubi \cite{Ka}. Roughly speaking, it is a specific choice of a Thom class of the stable normal bundle of each morphism of certain complex manifolds.

\begin{theorem}\label{3-3_analytic-Riemann-roch}
There exists an assignment $f\mapsto f_*$ sending a morphism $f\colon X\rightarrow Z$ of smooth projective varieties to a group homomorphism $f_*\colon K^0(X)\rightarrow K^0(Z)$ compatible with the natural homomorphism $K^0_{a}(-)\rightarrow K^0(-)$ and such that
$$\ch(f_*(\kappa))\cup\td(Z)=f_*(\ch(\kappa)\cup\td(X))$$ 
for each $\kappa\in K^0(X)$. 
	
\begin{proof}
One can refer to \cite[3.1]{AH61} and \cite[Theorem V.4.17]{Ka}. 
\end{proof}
\end{theorem}

The following statement is claimed in \cite[§4]{AH62} and known to experts. One can find proofs for (A) and (B) in \cite[Section IV.5]{Ka} and (C) is due to the specific choice of the Thom classes.

\begin{proposition}\label{3-2_push-forward-properties}
The assignment $f\mapsto f_*$ mentioned in Theorem \ref{3-3_analytic-Riemann-roch} also satisfies
	\begin{itemize}
		\item[(A)] $(f\circ g)_*=f_*\circ g_*$ for any $f\colon Z\rightarrow W,g\colon X\rightarrow Z$;
		\item[(B)] $f_*(\kappa_1\cup f^*\kappa_2)=f_*\kappa_1\cup\kappa_2$ for any $\kappa_1\in K^0(X),\kappa_2\in K^0(Z)$ and $f\colon X\rightarrow Z$;
		\item[(C)] $\pi_Z^*\circ f_*=(f\times\textup{id})_*\circ\pi_X^*$ for $f\colon X\rightarrow Z,\pi_X\colon X\times W\rightarrow X,\pi_Z\colon Z\times W\rightarrow Z$;
	\end{itemize}
where the maps involved are morphisms of smooth projective varieties.
\end{proposition}

\begin{proposition}\label{3-3_topo-K-FM-composition}
One has $\Phi_P^K\circ\Phi_Q^K=\Phi_{P\circ Q}^K$ for any $P\in D^b(Z\times W)$ and $Q\in D^b(X\times Z)$.
	
\begin{proof}
This is basically the same as \cite[Proposition 5.10]{HuFu} since the pushforward for topological $K$-theory satisfies all the necessary properties, cf.\ Proposition \ref{3-2_push-forward-properties}.
\end{proof}
\end{proposition}

\begin{proposition}\label{3-3_commutative}
Consider smooth projective varieties $X,Z$ and $P\in D^b(X\times Z)$, then there exists a commutative diagram
\begin{displaymath}\xymatrix{
			D^b(X)\ar[d]_{\Phi_P}\ar[rr]^{[-]}&&K^0(X)\ar[d]_{\Phi_P^K}\ar[rr]^{\mukai}&&H^{\textup{even}}(X,\bq)\ar[d]_{\Phi_P^H}\\
			D^b(Z)\ar[rr]^{[-]}&&K^0(Z)\ar[rr]^{\mukai}&&H^{\textup{even}}(Z,\bq)
}\end{displaymath}

\begin{proof}
Since the pushforward for topological $K$-theory satisfies the projective formula and the analytic Riemann--Roch theorem, one can prove it in the same way as \cite[Corollary 5.29]{HuFu}.
\end{proof}
\end{proposition}

\subsection{}
Orlov proved in \cite{Or97} that the kernel for a fully faithful Fourier--Mukai transform is unique. In general, the kernel is not unique but one always has

\begin{proposition}[{{\cite[Corollary 4.4]{CS12}}}]\label{3-4_class-of-kernel}
Consider smooth projective varieties $X,Z$ and $P,Q\in D^b(X\times Z)$ such that $\Phi_{P}\cong\Phi_{Q}$, then $\shh^i(P)\cong\shh^i(Q)$ for all $i\in\bz$. 
\end{proposition}

Proposition \ref{3-4_class-of-kernel} is usually applied together with the following criterion cf.\ \cite[Lemma 3.31]{HuFu} to show that the kernel of certain Fourier--Mukai transform is unique.

\begin{lemma}[Bridgeland]\label{3-4_point-sheaf-and-kernel}
Consider smooth projective varieties $X,Z$ and a Fourier--Mukai transform $\Phi_P\colon D^b(X)\rightarrow D^b(Z)$ such that for any point $x\in X$ there exists a coherent sheaf $\shk^x$ on $Z$ with $\Phi_{P}(k(x))\cong\shk^x[m]$ in $D^b(Z)$ for some $m\in\bz$, then $P\cong\shp[m]$ for some coherent sheaf $\shp$ on $X\times Z$ which is flat over $X$ via the canonical projection. 
\end{lemma}

\begin{example}\label{3-4_kernel-of-projection}
Let $Y$ be a cubic threefold, then the left mutation functor
$$\funct{L}_{\sho_Y(n)}\colon E\mapsto\textsf{C}(\ev\colon\bigoplus_m\Hom(\sho_Y(n),E[m])\otimes \sho_Y(n)[-m]\rightarrow E)$$
or more rigorously $j_*\circ\funct{L}_{\sho_Y(n)}$ is a Fourier--Mukai transform of kernel $(\sho(-n,n)\rightarrow\sho_{\Delta})$ where $\sho(a,b):=q^*\sho_X(a)\otimes p^*\sho_X(b)$ and $\sho_{\Delta}$ sits in degree $0$. Hence $j^*\cong\funct{L}_{\sho_Y}\circ\funct{L}_{\sho_Y(1)}$ or $j_*\circ j^*$ is also a Fourier--Mukai transform, with kernel
\begin{align*}
	U\cong(\sho(0,0)\rightarrow\sho_{\Delta})\circ(\sho(-1,1)\rightarrow\sho_{\Delta})\cong \left(\sho(-1,0)^{\oplus 5}\rightarrow\sho(-1,1)\oplus\sho(0,0)\rightarrow\sho_{\Delta}\right)
\end{align*}
where $\sho_{\Delta}$ sits in degree $0$. On the other hand, one has $j^*k(x)=\funct{L}_{\sho_Y}\circ\funct{L}_{\sho_Y(1)}k(x)\cong \shk^x[2]$ for a coherent sheaf $\shk^x$ on $Y$. This means that $U\cong\shu[2]$ in $D^b(Y\times Y)$ for some coherent sheaf $\shu$ on $Y\times Y$ by Lemma \ref{3-4_point-sheaf-and-kernel}. Hence $\Phi_P\cong\Phi_U$ implies that $P\cong U$ by Theorem \ref{3-4_class-of-kernel}.
\end{example}

\section{The group of Fourier--Mukai type autoequivalences}\label{SS_Main Theorem}
\subsection{}
Now we are prepared to work on $\Aut_{FM}(\scha_Y)$ for a cubic threefold $Y$. 

\begin{definition}
Given two cubic threefolds $Y_1$ and $Y_2$, an exact functor $\Phi\colon\scha_{Y_1}\rightarrow\scha_{Y_2}$ is said to be \defi{of Fourier--Mukai type of kernel $P$} if the composition
$$D^b(Y_1)\stackrel{j_1^*}{\longrightarrow}\scha_{Y_1}\stackrel{\Phi}{\longrightarrow}\scha_{Y_2}\stackrel{j_{2*}}{\longrightarrow}D^b(Y_2)$$
is isomorphic to the Fourier--Mukai transform $\Phi_P$, where $j_{2*}$ is the inclusion functor and $j_1^*$ is the left adjoint of the inclusion functor  $j_{1*}\colon \scha_{Y_1}\subset D^b(Y_1)$.
\end{definition}

Many natural functors in $\Aut(\scha_Y)$ appear to be Fourier--Mukai type functors.

\begin{example}
The identity $\funct{Id}_{\scha_Y}$ is a Fourier--Mukai type functor of kernel $U$ by Example \ref{3-4_kernel-of-projection}.
\end{example}

\begin{example}
The functors $\deriver f_*$ and $\derivel f^*$ are of Fourier--Mukai type for any $f\in\Aut(Y)$.
\end{example}

\begin{example}
The functors $\funct{T}$ and $\funct{S}$ are of Fourier--Mukai type according to \cite{HR19}.
\end{example}

Any composition of Fourier--Mukai type functors is of Fourier--Mukai type, so the set $\Aut_{FM}(\scha_Y)$ of (the isomorphic classes of) all Fourier--Mukai type autoequivalences of $\scha_Y$ has a natural monoid structure. Indeed, it is a group.

\begin{proposition}
The monoid $\Aut_{FM}(\scha_Y)$ is a group.

\begin{proof}
Given a Fourier--Mukai type autoequivalence $\Phi\in\Aut_{FM}(\scha_Y)$, it suffices to check that the composition $j_*\circ\Phi^{-1}\circ j^*\colon D^b(Y)\rightarrow D^b(Y)$ is a Fourier--Mukai transform.

Since $\funct{S}$ is a Fourier--Mukai type functor, the composition $\funct{S}\circ\Phi$ is of Fourier--Mukai type say with kernel $Q$. The right adjoint $\Phi_{Q_{\textrm{R}}}$ of $\Phi_Q$ is isomorphic to $j_*\circ\Phi^{-1}\circ\funct{S}^{-1}\circ j^!$. Set $\funct{S}_Y$ to be the Serre functor of $D^b(Y)$, then $j^!\cong \funct{S}\circ j^*\circ\funct{S}^{-1}_Y$. Thus $\Phi_{Q_{\textrm{R}}}\circ\funct{S}_Y\cong j_*\circ\Phi^{-1}\circ j^*$. It is known that $\funct{S}_Y$ is a Fourier--Mukai transform, so we are done.
\end{proof}
\end{proposition}

\subsection{}
Inspired by the work \cite{AT14}, we consider the group
$$K^0(\scha_Y):=\{\kappa\in K^0(Y)\,|\,\exists E\in\scha_Y\textup{ such that }\kappa=[E]\}$$
and fix the notation $\theta\colon K^0(\scha_Y)\subset K^0(Y)$ for convenience. 

\begin{lemma}\label{3-3_lemma}
Consider a cubic threefold $Y$, the functor $j_*\circ j^*\cong \Phi_U$ and its right adjoint $j_*\circ j^!\cong\Phi_{U_{\textup{\textrm{R}}}}$, then one has $\Phi^K_U\circ\theta=\Phi^K_{U_{\textup{\textrm{R}}}}\circ\theta=\theta$.

\begin{proof}
Given any $P\in D^b(Y\times Y)$, one has the following commutative diagram
\begin{displaymath}
	\xymatrix{
		\scha_Y\ar[rr]^{j_*}\ar[d]&& D^b(Y)\ar[rr]^{\Phi_P}\ar[d]&& D^b(Y)\ar[d]\\
		K^0(\scha_Y)\ar[rr]^{\theta}&& K^0(Y)\ar[rr]^{\Phi^K_P}&& K^0(Y)
	}
\end{displaymath}
where the left hand square is commutative by definition and the right hand square is commutative due to Proposition \ref{3-3_commutative}. Then our statement follows immediately from plugging $U$ and $U_{\textrm{R}}$.
\end{proof}
\end{lemma}

\begin{proposition}
Consider a cubic threefold $Y$ and a Fourier--Mukai type autoequivalence $\Phi$ on $\scha_Y$ of kernel $P$, then the topological $K$-theoretic map $\Phi_P^K\colon K^0(Y)\rightarrow K^0(Y)$ restricts to a group automorphism $\Phi^K\colon K^0(\scha_Y)\rightarrow K^0(\scha_Y)$.
	
\begin{proof}
The right adjoint functor of $\Phi_P\cong j_*\circ\Phi\circ j^*$ is $\Phi_{P_{\textrm{R}}}\cong j_*\circ\Phi^{-1}\circ j^!$, so one has
$$\Phi_P\circ\Phi_{P_{\textrm{R}}}\cong j_*\circ j^!\cong\Phi_{U_{\textrm{R}}}\quad\textup{and}\quad\Phi_{P_\textrm{R}}\circ\Phi_P\cong j_*\circ j^*\cong\Phi_U$$ 
It follows form Proposition \ref{3-3_topo-K-FM-composition} that $\Phi^K_P\circ\Phi^K_{P_\textrm{R}}=\Phi^K_{U_{\textrm{R}}}$ and $\Phi^K_{P_\textrm{R}}\circ\Phi^K_P=\Phi^K_U$. Then one can see 
\begin{center}
$\Phi^K_P\circ\theta\circ\Phi^K_{P_\textrm{R}}\circ\theta=\Phi^K_P\circ\Phi^K_{P_\textrm{R}}\circ\theta=\theta\quad$ and $\quad\Phi^K_{P_\textrm{R}}\circ\theta\circ\Phi^K_P\circ\theta=\Phi^K_{P_\textrm{R}}\circ\Phi^K_P\circ\theta=\theta$
\end{center}
due to Lemma \ref{3-3_lemma}. So the homomorphism $\Phi_P^K\circ\theta\colon K^0(\scha_Y)\hookrightarrow K^0(Y)$ injectively maps onto the subgroup $K^0(\scha_Y)$. In other words, $\Phi_P^K$ restricts to a group automorphism $\Phi^K$ of $K^0(\scha_Y)$.
\end{proof}
\end{proposition}

Given $\Phi_1,\Phi_2\in\Aut(\scha_Y)$ with kernel $P_1$ and $P_2$ respectively, one has 
$$\Phi_{P_1}^K\circ\Phi_{P_2}^K=\Phi_P^K$$ 
for any $P\in D^b(Y\times Y)$ satisfying $\Phi_P\cong\Phi_{P_1}\circ\Phi_{P_2}$ due to Proposition \ref{3-4_class-of-kernel}. So their restrictions on $K^0(\scha_Y)$ are the same i.e.\ $\Phi_1^K\circ\Phi_2^K=(\Phi_1\circ\Phi_2)^K$. This ensures a group homomorphism $$\rho\colon\Aut_{FM}(\scha_Y)\rightarrow\Aut(K^0(\scha_Y)),\quad\Phi\mapsto\Phi^K$$
where $\Aut(K^0(\scha_Y))$ is the group of all group automorphisms of $K^0(\scha_Y)$.

\subsection{}
To understand $\rho$, we need to know $K^0(\scha_Y)$ better. Choose any base point on $Y$, the analytic Riemann--Roch Theorem \ref{3-3_analytic-Riemann-roch} gives an Euler pairing on $K^0(\scha_Y)$ which is compatible with the usual Euler pairing on $K^0_{a}(\scha_Y)$. These Euler pairings will be both denoted by $\chi$.

\begin{proposition}
Consider any line $\ell\subset Y$, then there exists an isomorphism $$K^0(\scha_Y)\cong\bz[\shi_{\ell}]\oplus\bz[\funct{S}\shi_{\ell}]$$ 
and the Euler pairing $\chi(-,-)$ on $K^0(\scha_Y)$ has the form 
$$\begin{pmatrix}
		-1&-1\\0&-1
\end{pmatrix}$$
with respect to this basis.
	
\begin{proof}
As the integral cohomology of $Y$ are torsion-free, the topological $K$-theory of $Y$ is torsion-free as well using \cite[2.5]{AH61}. Then one has $N(\scha_Y)=\textup{im}(K^0_{a}(\scha_Y)\rightarrow K^0(\scha_Y))=K^0(\scha_Y)$ as cubic threefolds also satisfy the integral Hodge conjecture. The identification is natural and respects the Euler pairings, so one concludes by Proposition \ref{2_numerical-Grothendieck-group}.
\end{proof}
\end{proposition}

Since $\funct{S}^3=[5]$ and $\chi(\shi_{\ell},\shi_{\ell})=-1$, the proposition implies that
$$\{\kappa\in K^0(\scha_Y)\,|\,\chi(\kappa,\kappa)=-1\}=\{[\funct{S}^r\shi_{\ell}]\,|\,r=0,1,2,3,4,5\}$$
An automorphism $\Phi^K$ lying in the image of $\rho\colon\Aut_{FM}(\scha_Y)\rightarrow\Aut(K^0(\scha_Y))$ is determined by the image $\Phi^K([\shi_{\ell}])$. Indeed, one has
$$\Phi^K[\funct{S}\shi_{\ell}]=\Phi^K\funct{S}^K[\shi_{\ell}]=\funct{S}^K\Phi^K[\shi_{\ell}]$$
as $\Phi\circ\funct{S}=\funct{S}\circ\Phi$ for any autoequivalence $\Phi$ on $\scha_Y$. In particular, one has $\Phi^K([\shi_{\ell}])=[\funct{S}^r\shi_{\ell}]$ for some $r=0,1,2,3,4,5$ and thus $\image(\rho)\cong\bz/6\bz$. To identify $\ker(\rho)$, one observes

\begin{proposition}
Consider a Fourier--Mukai type autoequivalence $\Phi\colon\scha_Y\rightarrow\scha_Y$, then 
$$\Phi\cong\funct{S}^r\circ[2m]\circ\Psi$$ 
for a unique integer $m$, a unique $r\in\{0,1,2,3,4,5\}$ and a unique Fourier--Mukai type functor $\Psi\colon\scha_Y\rightarrow\scha_Y$ sending the ideal sheaf of a line to the ideal sheaf of some line.
	
\begin{proof}
The existence of the functor $\Psi$ is guaranteed by Proposition \ref{2-3_decomposition}. The uniqueness for all of $m,r,\Psi$ follows from the fact that $\image(\rho)\cong\bz/6\bz$ and $\funct{S}^3\cong[5]$.
\end{proof}
\end{proposition}

This proposition implies a bijective group homomorphism
$$\Aut_{FM}(\scha_Y)\rightarrow \bz\times G,\quad\Phi=\funct{S}^r\circ[2m]\circ\Psi\mapsto(5r+6m,\Psi)$$
where 
$$G=\{\Psi\in\Aut_{FM}(\scha_Y)\,|\,\Psi(\shi_{\ell})\cong\shi_{\ell'}\textup{ for any line }\ell\textup{ and some line }\ell'\}$$		
One can associate $\Psi\in G$ with an element in $\Aut(F(Y))\cong\Aut(Y)$ by sending $\ell$ to $\ell'$, where $\ell'$ is the unique line such that $\shi_{\ell'}=\Psi(\shi_{\ell})$. This gives a group homomorphism $\rho_0\colon G\rightarrow\Aut(Y)$ which is only known to be surjective so far by $f\mapsto f^*$. The next step is to show that $\rho_0$ is an isomorphism, and it suffices to show that
$$G_0=\{\Psi\in\Aut_{FM}(\scha_Y)\,|\,\Psi(\shi_{\ell})\cong\shi_{\ell}\textup{ for any line }\ell\subset Y\}$$	 
is the trivial subgroup of $\Aut_{FM}(\scha_Y)$. 

\subsection{}
Proposition \ref{2-3_point-identification} allows us to make the following observation.

\begin{proposition}
	One has $\Psi(\shk^x)\in\{\shk^{y}[2n]\,|\,y\in Y,n\in\bz\}$ for any $\Psi\in G_0$.
\end{proposition}

Since $\shk^x[2]\cong j^*k(x)$, one sees $\Hom(\shk^x,\shi_{\ell}[n])\cong\Hom(k(x),\shi_{\ell}[n+2])$ by $j^* \dashv j_*$. Hence one can detect whether $x\in\ell$ by the vector spaces $\Hom(\shk^x,\shi_{\ell}[n])$.

\begin{proposition}
Any autoequivalence $\Psi\in G_0$ satisfies $\Psi(\shk^x)\cong\shk^x$ for any point $x\in Y$.
	
\begin{proof}
Given a line $\ell\subset Y$ and a point $x\in Y$, it is standard to compute that
$$x\notin\ell\Rightarrow \Hom(k(x),\shi_{\ell}[n+2])\cong\left\{\begin{aligned}
	\bc&&n=1\\
	0&&n\neq 1	
\end{aligned}\right.$$
and
$$x\in\ell\Rightarrow \Hom(k(x),\shi_{\ell}[n+2])\cong\left\{\begin{aligned}
	\bc\,\,\,\,\,\, &&&n=0\\
	\bc^{\oplus 2} &&&n=1\\
	0\,\,\,\,\,\,  &&&n\neq 0,1
\end{aligned}\right.$$
Therefore, one has $\Psi(\shk^x)\cong\shk^{y_x}$ for some point $y_x$ such that $y_x\in\ell\Leftrightarrow x\in\ell$ for any $\ell\subset Y$. Since a point on $Y$ is fixed by lines, one must have $x=y_x$.
\end{proof}
\end{proposition}

According to Lemma \ref{3-4_point-sheaf-and-kernel}, the kernel of $\Phi\in G_0$ is isomorphic to $\shp[2]$ for a coherent sheaf $\shp$ on $Y\times Y$ flat over the first factor and $\iota_x^*\shp\cong\shk^x$ for each $x\in Y$ where $\iota_x$ denotes the map $Y\rightarrow Y\times Y,y\mapsto(x,y)$. Recall that the moduli space $M_Y(v)$ is a fine moduli space, so there is an isomorphism of functors $\alpha\colon\mathfrak{M}\rightarrow\Hom(-,M_Y(v))$ where $\mathfrak{M}$ is the moduli functor associated with $M_Y(v)$. In particular, one has 
$$\alpha_Y\colon\mathfrak{M}(Y)=\Hom(Y,M_Y(v))$$
and can see that the morphism $g=\alpha_Y(\shp)$ satisfies $g(x)=x$ for any point $x\in Y$. Since both $Y$ and $M_Y(v)$ are smooth projective varieties, one has $g^{-1}(x)\cong\spec(k(x))$ and that $g$ is projective. So $g$ is just the closed embedding $Y\hookrightarrow M_Y(v)$ cf.\ \cite[18.12.6]{GroEGAIV4}.

As in Example \ref{3-4_kernel-of-projection}, we denote the sheaf associated with $\funct{Id}_{\scha_Y}$ by $\shu$. Then $\alpha_Y(\shp)=\alpha_Y(\shu)$ indicates that the two flat families $\shp$ and $\shu$ only differ by a line bundle $\sho_Y(m)$ on $Y$.

\begin{proposition}
Any $\Psi\in G_0$ is isomorphic to $\funct{Id}_{\scha_Y}$.
	
\begin{proof}
The kernel of $\Psi$ is of the form $\shp[2]$ for a flat family $\shp$ of stable sheaves such that $\iota_x^*\shp\cong\shk^x$ for each $x\in Y$. Then there exists some $m\in\bz$ such that $\shp\cong\shu\otimes q^*\sho_Y(m)$.
		
It follows that  
$$\shk^x\cong \Phi_{\shp[2]}\shk^x\cong p_*(\shu[2]\otimes q^*(\shk^x\otimes\sho_Y(m)))\cong j_*j^*(\shk^x\otimes\sho_Y(m)) \cong\funct{T}^m(\shk^x)$$ 
One can see $m=0$ by dividing into three cases:
\begin{itemize}
	\item $m=3k$ implies that $\shk^x\cong\funct{T}^{3k}\shk^x\cong\shk^x[2k]$. So $k=0$.
	\item $m=3k+1$ implies that $\shk^x\cong\funct{T}^{3k+1}\shk^x\cong\funct{T}\shk^x[2k]$, which is impossible.
	\item $m=3k+2$ implies that $\shk^x\cong\funct{T}^{3k+2}\shk^x\cong\funct{T}^2\shk^x[2k]$, which is impossible.
\end{itemize}
Therefore, one can only have $m=0$ and $\Phi_{\shp[2]}\cong\Phi_{\shu[2]}\cong j_*\circ j^*$. It means that $\Psi\cong\funct{Id}_{\scha_Y}$. 
\end{proof}
\end{proposition}

In conclusion, we have obtained the following theorem.

\begin{theorem}\label{5_main-theorem}
Consider a cubic threefold $Y$, then there exists a group isomorphism $$\Aut_{FM}(\scha_Y)\cong\bz\times\Aut(Y),\quad\Phi\mapsto(5r+6m,f)$$
where $\Phi=\funct{S}^r\circ[2m]\circ f^*$ for some $m\in\bz$, $r\in\{0,1,2,3,4,5\}$ and $f\in\Aut(Y)$.
\end{theorem}

\section{The Fourier--Mukai type categorical Torelli theorem revisited}\label{SS_FM categorical Torelli}

\subsection{} Here we present our proof for the Fourier--Mukai version of categorical Torelli theorem.

\begin{theorem}\label{5_another-proof}
Consider two cubic threefolds $Y_1$ and $Y_2$, and assume that there exists a Fourier--Mukai type equivalence $\scha_{Y_2}\cong\scha_{Y_2}$, then $Y_1\cong Y_2$.
\end{theorem}

Given such an equivalence $\Phi$ and name the associated Fourier--Mukai transform $\Phi_P$, we will have a cohomological Fourier--Mukai map
$$\Phi_P^H\colon H^{\textup{odd}}(Y_1,\bq)\rightarrow H^{\textup{odd}}(Y_2,\bq)$$
Since the only non-trivial odd-degree cohomology group of a cubic threefold $Y$ is $H^3$, it then restricts to a group homomorphism $\Phi_P^H\colon H^3(Y_1,\bq)\rightarrow H^3(Y_2,\bq)$.

We will show that the group homomorphism is indeed a Hodge isomorphism, so that we can apply the arguments in \cite{OR18} to see that $\Phi_P^H$ restricts again to a Hodge isometry 
$$\Phi_P^H\colon H^3(Y_1,\bz)\rightarrow H^3(Y_2,\bz)$$
In this case, one can conclude by the classic Torelli theorem for cubic threefolds.

\subsection{}
The key point is the following lemma, which is pointed out to me by Daniel Huybrechts.

\begin{lemma}\label{5_zero-homomorphism}
Consider a cubic threefold $Y$ and any two coherent sheaves $\she$ and $\shf$ on $Y$, then the cohomological Fourier--Mukai map $\Phi^H_{\she\boxtimes\shf}\colon H^3(Y,\bq)\rightarrow H^3(Y,\bq)$ is the zero map. Here the kernel $\she\boxtimes\shf:=q^*\she\otimes p^*\shf$ as before.

\begin{proof}
Since Betti cohomology is a classic Weil cohomology theory, one has
$$\Phi^H_{\she\boxtimes\shf}(v)=p_*(q^*(v\cup \mukai(\she))\cup p^*\mukai(\shf))=\left(\int_Yv\cup \mukai(\she)\right)\cup\mukai(\shf)$$ 
for any class $v\in H^3(Y,\bq)$. The cohomology class $v\cup\mukai(\she)$ has no components in the top cohomology group, so it is zero after the trace map. Hence $\Phi^H_{\she\boxtimes\shf}(v)$ has to be zero.
\end{proof}
\end{lemma}

Since $p_*,q^*,\cup$ are additive, one can see that 
$$
p_*(q^*v\cup(\mukai([E]+[F])))=p_*(q^*v\cup(\mukai([E])))+ p_*(q^*v\cup(\mukai([F])))
$$
for any $E,F\in D^b(Y\times Y)$ and cohomology class $v$. Moreover, one has the following observations:

\begin{corollary}\label{5_corollary-1}
One has $\Phi^H_U=\Phi^H_{\sho_{\Delta}}=\textup{Id}_{H^3}$ for the Fourier--Mukai transform $j_*\circ j^*\cong\Phi_U$. 
	
\begin{proof}
According to Example \ref{3-4_kernel-of-projection}, one has $$[U]=[\sho_{\Delta}]+5[\sho(-1,0)]-[\sho(-1,1)]-[\sho(0,0)]$$ 
in the topological $K$-group $K^0(Y\times Y)$, which implies $\Phi_U^H=\Phi^H_{\sho_{\Delta}}=\textup{Id}_{H^3}$ by Lemma \ref{5_zero-homomorphism} and the fact that the Mukai vector $\mukai\colon K^0(Y\times Y)\rightarrow H^{\textup{even}}(Y\times Y)$ is also additive.
\end{proof}
\end{corollary}

The functor $j_*\circ j^!$ is the right adjoint of the Fourier--Mukai transform $\Phi_U\cong j_*\circ j^*$, so one has $j_*\circ j^!\cong\Phi_{U_{\textup{\textrm{R}}}}$ for the kernel $U_{\textup{\textrm{R}}}=U^\vee\otimes q^*\sho_Y(-2)[3]$.

\begin{corollary}\label{5_corollary-2}
One has $\Phi^H_{U_{\textup{\textrm{R}}}}=\Phi^H_{\sho_{\Delta}}=\textup{Id}_{H^3}$ for the Fourier--Mukai transform $j_*\circ j^!\cong\Phi_{U_{\textup{\textrm{R}}}}$. 
	
\begin{proof}
The derived dual of a locally free sheaf on $Y$ is just the usual dual, so one has 
$$\sho(a,b)^\vee\otimes q^*\sho_Y(-2)[3]\cong\sho(-a-2,-b)[3]$$ 
Since $\Phi_{\sho_{\Delta}}\cong\funct{Id}_{D^b(Y)}$, its right adjoint functor is also the identity functor. Hence one has
$$\sho_{\Delta,\textup{\textrm{R}}}:=\sho^\vee_{\Delta}\otimes q^*\sho_Y(-2)[3]\cong\sho_{\Delta}$$ 
due to Proposition \ref{3-4_class-of-kernel}. Therefore, one can conclude by Lemma \ref{5_zero-homomorphism} and the isomorphism 
$$U\cong\left(\sho(-1,0)^{\oplus 5}\rightarrow\sho(-1,1)\oplus\sho(0,0)\rightarrow\sho_{\Delta}\right)$$
that one has $\Phi^H_{U_{\textup{\textrm{R}}}}=\Phi^H_{\sho_{\Delta,\textup{\textrm{R}}}}=\Phi_{\sho_{\Delta}}^H=\textup{Id}_{H^3}$.
\end{proof}
\end{corollary}

\begin{proposition}
Consider two cubic threefolds $Y_1,Y_2$ and a Fourier--Mukai type equivalence $\scha_{Y_1}\cong\scha_{Y_2}$ of kernel $P$, then $\Phi_P^H\colon H^3(Y_1,\bq)\rightarrow H^3(Y_2,\bq)$ is a group isomorphism.
	
\begin{proof}
One has $\Phi_{P_\textup{\textrm{R}}\circ P}\cong\Phi_{P_\textup{\textrm{R}}}\circ\Phi_P\cong j_{1*}\circ j^*_1$, so by Corollary \ref{5_corollary-1} the map 
$$\Phi^H_{P_\textup{\textrm{R}}}\circ\Phi^H_P\cong\Phi_{P_\textup{\textrm{R}}\circ P}^H\colon H^3(Y_1,\bq)\rightarrow H^3(Y_1,\bq)$$
is the identity map. Similarly, the map $\Phi^H_{P}\circ\Phi^H_{P_\textup{\textrm{R}}}$ is also the identity due to Corollary \ref{5_corollary-2}.
\end{proof}
\end{proposition}

In this case, the cohomological Fourier--Mukai map $\Phi_P^H$ restricts to isomorphisms
$$\bigoplus_{p-q=-1}H^{p,q}(Y_1)\cong \bigoplus_{p-q=-1}H^{p,q}(Y_2)\quad\textup{and}\quad \bigoplus_{p-q=1}H^{p,q}(Y_1)\cong \bigoplus_{p-q=1}H^{p,q}(Y_2)$$
which implies that $\Phi_P^H\colon H^3(Y_1,\bq)\rightarrow H^3(Y_2,\bq)$ is a Hodge isomorphism. The remaining part of our proof is the same as that for \cite[Proposition 2.1]{OR18}.

\end{document}